\documentclass[11pt]{amsart}
\usepackage{graphicx}
\usepackage{amsmath}
\usepackage{amssymb}
\usepackage{amsfonts}
\usepackage{tikz}
\usetikzlibrary{tqft,shapes}
\usepackage{cite}
\usepackage{tikz-cd}
\usepackage[margin=1.40in]{geometry}

\newtheorem{theorem}{Theorem}[subsection]

\newtheorem{lemma}[theorem]{Lemma}

\newtheorem{remark}[theorem]{Remark}

\theoremstyle{definition}
\newtheorem{definition}[theorem]{Definition}
\theoremstyle{remark}
\numberwithin{equation}{subsection}

\begin{document}

\title{On twisted torsion of compact $3$-manifolds}
\author{ESMA D\.IR\.ICAN ERDAL}

\address[1]{%
Deptartment of Mathematics\\
I\c{s}{\i}k University\\
 \.{I}stanbul\\
  T\"{u}rkiye}
\email{esma.diricanerdal@isikun.edu.tr}

\subjclass[2010]{Primary 55U99; Secondary (optional) 18G99, 57Q10}
\keywords{Reidemeister torsion, twisted invariants, connected sum, 3-manifolds.}

\begin{abstract}
In this paper, we examine the behavior of twisted Reidemeister torsion for compact, connected, oriented $3$-manifolds under the connected sum operation. Under appropriate assumptions, we establish a multiplicative relation for the torsion associated with representations into a complex reductive algebraic group. Our results are consistent with the known behavior of torsion-type invariants and provide a formulation that may be useful for further investigations of decomposition properties in $3$-manifold topology.
\end{abstract}
\maketitle

\section{Introduction}
\label{Sec:1}

Reidemeister torsion is a classical invariant in topology, originally introduced by Reidemeister and Franz in the study of lens spaces \cite{Reidemeister-Franz,Franz}. It was later shown to be a topological invariant in broader contexts: Kirby and Siebenmann established invariance for piecewise-linear manifolds \cite{RCLC}, and Chapman extended this to general simplicial complexes \cite{Chap}. These developments ensured that the classification of lens spaces via Reidemeister torsion is topologically well-defined.

Beyond its foundational role, Reidemeister torsion has found applications in several areas of topology. Milnor employed it in the study of triangulations, constructing homeomorphic but combinatorially distinct simplicial complexes and thereby disproving the Hauptvermutung \cite{Milnor2,Milnor}. Moreover, its relationship with the Alexander polynomial has made it a useful tool in knot theory, and variants of torsion appear in the construction of more refined invariants such as the Casson--Walker--Lescop invariant \cite{Lescop}.

In this work we consider compact, connected, oriented $3$-manifolds, possibly with non-empty boundary, equipped with a fixed CW-decomposition. Let $G$ be a complex reductive algebraic group (e.g. $SL_n(\mathbb{C})$ or $PSL_n(\mathbb{C})$) with Lie algebra $\mathfrak{g}$, and let $\rho:\pi_1(X)\to G$ be a representation. Associated to these data one can define twisted homology groups and the corresponding twisted Reidemeister torsion $\mathbb{T}_\rho$.

Understanding how $\mathbb{T}_\rho$ behaves under the connected sum operation is a natural and important question, since it clarifies how a global torsion invariant reflects a manifold's decomposition into simpler pieces. In dimension two, Heusener and Porti showed that the twisted Reidemeister torsion of a four-holed sphere can be computed by viewing it as the union of two pairs of pants glued along two points \cite{Porti2}. In dimension three, for closed manifolds, Porti and Yoon computed the twisted Reidemeister torsion of the connected sum of knot exteriors \cite{Porti3}. An explicit, concrete multiplicativity formula for the twisted Reidemeister torsion of a connected sum of compact $3$-manifolds with boundary appears, however, to be missing from the literature. The purpose of the present paper is to fill this gap.

Our strategy relies on a classical fact due to Gross \cite{JLGROSS}: 
\emph{if $M$ is a compact $3$-manifold with connected non-vacuous boundary and $N$ is a $3$-manifold with vacuous boundary, and if $N^{\ast}$ denotes $N$ with the interior of a $3$-cell removed, then the connected sum $M\# N$ is homeomorphic to the \emph{disk sum} $M\bigtriangleup N^{\ast},$ namely $M\# N\approx M\bigtriangleup N^{\ast}$. }
Combining this homeomorphism with the multiplicativity formula for the disk sum established in \cite[Theorem 1.2]{EDErdal}, we obtain the following theorem, which is the main result of this paper.

\begin{theorem}\label{primetheo}
Let $M$ be a compact $3$-manifold with connected non-vacuous boundary which is not spherical, and let $N$ be another compact $3$-manifold with vacuous boundary. Let $N^{\ast}$ denote $N$ with the interior of a $3$-cell removed, and let
$$\varrho:\pi_1(M\# N)\longrightarrow G$$
be a given representation. Set $\sigma=[\varphi^{-1}]\circ\varrho:\pi_1(M\bigtriangleup N^{\ast})\to G$, where $\varphi$ is the homeomorphism 
$M\# N\approx M\bigtriangleup N^{\ast}$, and let $\psi_M=\sigma_{|_{\pi_1(M)}}$ and $\psi_{N^{\ast}}=\sigma_{|_{\pi_1(N^{\ast})}}$ be the restrictions of 
$\sigma.$ For a given homology basis $\{\mathbf{h}_i^{M\# N}\}_{i=0}^3$ of $H_i(M\# N;\mathfrak{g}_{\mathrm{Ad}_{\varrho}})$, there exist homology bases $\{\mathbf{h}_i^{M}\}_{i=0}^3$ and $\{\mathbf{h}_i^{N^{\ast}}\}_{i=0}^3$ of $H_i(M;\mathfrak{g}_{\mathrm{Ad}_{\psi_M}})$ and $H_i(N^{\ast};\mathfrak{g}_{\mathrm{Ad}_{\psi_{N^{\ast}}}})$ respectively, such that
\begin{equation}\label{mainformula}
\mathbb{T}_\varrho\big(M\# N,\{\mathbf{h}_i^{M\# N}\}_{i=0}^3\big)
=\mathbb{T}_{\psi_M}\big(M,\{\mathbf{h}_i^{M}\}_{i=0}^3\big)\;
\mathbb{T}_{\psi_{N^{\ast}}}\big(N^{\ast},\{\mathbf{h}_i^{N^{\ast}}\}_{i=0}^3\big).
\end{equation}
\end{theorem}

\begin{remark}
Formula \eqref{mainformula} is a pure product, with \emph{no corrective term}. This may seem surprising in view of the general Mayer--Vietoris multiplicativity formula of Theorem \ref{prt1}(iii), where a correction term $\mathbb{T}(\mathcal{H}_\ast,\{\mathbf{h}_{\ast\ast}\})$ typically appears. The reason it disappears here is that, once the connected sum is rewritten as a disk sum, the two pieces $M$ and $N^{\ast}$ are glued along a single topological disk $\mathbb{D}^2$ rather than along a sphere. Since $\mathbb{D}^2$ is contractible, its twisted torsion (computed with respect to the natural basis coming from a point) equals $1$ by \cite[Lemma 3.1]{EDErdal}, and the Mayer--Vietoris correction term coming from the disk turns out to be trivial as well. This is precisely the mechanism that makes Theorem \ref{primetheo} an effective computational tool.
\end{remark}

The behavior of Reidemeister torsion under topological decompositions has been studied in several contexts. In the classical (untwisted) setting, foundational results can be found in the work of Milnor, who established the connection with Whitehead torsion and its applications to triangulations. A systematic treatment, including refinements addressing indeterminacy, was later developed by Turaev, who introduced Euler structures and homology orientations to obtain more canonical formulations. In the twisted setting associated with representations of the fundamental group, torsion has been studied extensively for $3$-manifolds and their geometric structures; the work of Porti and others has highlighted the role of twisted torsion in relation to representation varieties and geometric invariants, typically formulated via Mayer--Vietoris sequences at the level of chain complexes. Our Theorem \ref{primetheo} fits into this general framework, adapted specifically to the connected sum decomposition of compact, oriented $3$-manifolds.

\subsection{Structure of the paper}
The paper is organized as follows. In Section~\ref{Sec:1} (this section) we introduce the problem and state the main result. In Section~\ref{Sec:2} we recall the necessary background on twisted Reidemeister torsion, most of which is classical and collected here for the reader's convenience and we also recall the disk sum operation and Theorem 1.2 of \cite{EDErdal}. In Section~\ref{Sec:3} 
we prove Theorem \ref{primetheo} by using the results given in Section~\ref{Sec:2}.

\subsection{Notation and conventions}
Throughout the paper, every $3$-manifold is triangulated, oriented, connected and compact. Unless stated otherwise, $M$ denotes a $3$-manifold with connected non-vacuous boundary which is not spherical, and $N$ denotes a $3$-manifold with vacuous (empty) boundary; $N^{\ast}$ is $N$ with the interior of a $3$-cell removed. The group $G$ is a complex reductive algebraic group, typically $SL_n(\mathbb{C})$ or $PSL_n(\mathbb{C})$, with Lie algebra $\mathfrak{g}$. For a finite CW-complex $X$ with cell-decomposition $K$ and a representation $\rho:\pi_1(X)\to G$, we write $C_{\ast}(K;\mathfrak{g}_{\mathrm{Ad}_\rho})$ for the associated twisted chain complex, and
$$H_{p}(K;\mathfrak{g}_{\mathrm{Ad}_\rho})=Z_{p}(K;\mathfrak{g}_{\mathrm{Ad}_\rho})/B_{p}(K;\mathfrak{g}_{\mathrm{Ad}_\rho})$$
for the $p$-th twisted homology group, where
$$Z_{p}(K;\mathfrak{g}_{\mathrm{Ad}_\rho})=\mathrm{Ker}\{\partial_{p}\otimes\mathrm{id}\},\qquad
B_{p}(K;\mathfrak{g}_{\mathrm{Ad}_\rho})=\mathrm{Im}\{\partial_{p+1}\otimes\mathrm{id}\}.$$
Given a homology basis $\{\mathbf{h}^{M}_p\}_{p=0}^3$ of $H_p(M;\mathfrak{g}_{\mathrm{Ad}_{\rho}})$, we write $\mathbb{T}_{\rho}(M,\{\mathbf{h}^{M}_p\}_{p=0}^3)$ for the twisted Reidemeister torsion of $M$ with respect to $\rho$; the precise definition is recalled in Section~\ref{Sec:2}.
\section{Preliminaries}\label{Sec:2}
\subsection{The Twisted Reidemeister Torsion}
\label{Subsec:2.1}

Consider an $n$-dimensional CW-complex $X$ with universal covering $\widetilde{X}$. Let $\mathcal{B}$ denote the non-degenerate Killing form on $\mathfrak{g}$, given by $\mathcal{B}(A,B)=4\cdot\mathrm{Trace}(AB)$. For a representation $\rho:\pi_1(X)\to G$, we let $\pi_1(X)$ act on $\mathfrak{g}$ via the adjoint of $\rho$.

Let $K$ be a cell-decomposition of $X$ and $\widetilde{K}$ a lift of $K$. Using the cellular chain complex $C_{\ast}(\widetilde{K};\mathbb{Z})$, the twisted chain complex is defined by
\begin{equation}\label{defchn1}
C_{\ast}(K;\mathfrak{g}_{\mathrm{Ad}_\rho}):=C_{\ast}(\widetilde{K};\mathbb{Z})\otimes\mathfrak{g}\big/\!\sim,
\end{equation}
where $\sigma\otimes t\sim \gamma\cdot\sigma\otimes\gamma\cdot t$ for every $\gamma\in\pi_1(X)$; here $\pi_1(X)$ acts on $\widetilde{X}$ by deck transformations, and on $\mathfrak{g}$ by the adjoint action.

Let $\{e_1^p,\ldots,e_{m_p}^p\}$ be the $p$-cells of $K$, forming a $\mathbb{Z}$-basis of $C_p(K;\mathbb{Z})$. Choosing a lift $\widetilde{e}^p_j$ of each $e^p_j$ gives a $\mathbb{Z}[\pi_1(X)]$-basis $c_p=\{\widetilde{e^p_j}\}_{j=1}^{m_p}$ of $C_p(\widetilde{K};\mathbb{Z})$. Fixing a $\mathcal{B}$-orthonormal basis $\mathcal{A}=\{\mathfrak{a}_k\}_{k=1}^{\dim\mathfrak{g}}$ of $\mathfrak{g}$, we obtain a \emph{geometric basis} $\mathbf{c}_p=c_p\otimes_{\rho}\mathcal{A}$ of $C_p(K;\mathfrak{g}_{\mathrm{Ad}_\rho})$.

Consider the chain complex $C_{\ast}:=C_{\ast}(K;\mathfrak{g}_{\mathrm{Ad}_\rho})$,
\begin{equation*}
0\to C_{n}(K;\mathfrak{g}_{\mathrm{Ad}_\rho})\to C_{n-1}(K;\mathfrak{g}_{\mathrm{Ad}_\rho})\to\cdots\to C_{0}(K;\mathfrak{g}_{\mathrm{Ad}_\rho})\to 0.
\end{equation*}
For $p\in\{0,\ldots,n\}$, let $H_p(C_\ast)=Z_p(C_\ast)/B_p(C_\ast)$ be the $p$-th homology, viewed as a $\mathbb{C}$-vector space, where
$$B_p(C_\ast)=\mathrm{Im}\{\partial_{p+1}:C_{p+1}\to C_p\},\qquad
Z_p(C_\ast)=\mathrm{Ker}\{\partial_p:C_p\to C_{p-1}\}.$$
This gives the short exact sequences
\begin{equation}\label{Equation1}
0\longrightarrow Z_p(C_\ast)\stackrel{\imath}{\hookrightarrow} C_p(C_\ast)\stackrel{\partial_p}{\longrightarrow} B_{p-1}(C_\ast)\longrightarrow 0,
\end{equation}
\begin{equation}\label{Equation2}
0\longrightarrow B_p(C_\ast)\stackrel{\imath}{\hookrightarrow} Z_p(C_\ast)\stackrel{\varphi_p}{\longrightarrow} H_p(C_\ast)\longrightarrow 0,
\end{equation}
where $\imath$ and $\varphi_p$ denote inclusion and natural projection, respectively.

Choosing sections $s_p:B_{p-1}(C_\ast)\to C_p(C_\ast)$ of $\partial_p$ and $\ell_p:H_p(C_\ast)\to Z_p(C_\ast)$ of $\varphi_p$, the Splitting Lemma applied to \eqref{Equation1}--\eqref{Equation2} gives
\begin{equation}\label{Equation0}
C_p(C_\ast)=B_p(C_\ast)\oplus\ell_p(H_p(C_\ast))\oplus s_p(B_{p-1}(C_\ast)).
\end{equation}
Thus, if $\mathbf{b}_p$ and $\mathbf{h}_p$ are bases of $B_p(C_\ast)$ and $H_p(C_\ast)$, the disjoint union $\mathbf{b}_p\sqcup\ell_p(\mathbf{h}_p)\sqcup s_p(\mathbf{b}_{p-1})$ is a new basis of $C_p(C_\ast)$, and we may compare it to the geometric basis $\mathbf{c}_p$.

\begin{definition}
Let $[\mathbf{b}_p\sqcup\ell_p(\mathbf{h}_p)\sqcup s_p(\mathbf{b}_{p-1}),\mathbf{c}_p]$ denote the determinant of the transition matrix between the two bases of $C_p(C_\ast)$ above. The \emph{twisted Reidemeister torsion} of $C_{\ast}(K;\mathfrak{g}_{\mathrm{Ad}_\rho})$ is
$$\mathbb{T}(C_{\ast}(K;\mathfrak{g}_{\mathrm{Ad}_\rho}),\{\mathbf{c}_p\}_{p=0}^{n},\{\mathbf{h}_p\}_{p=0}^{n})
=\prod_{p=0}^n\big[\mathbf{b}_p\sqcup\ell_p(\mathbf{h}_p)\sqcup s_p(\mathbf{b}_{p-1}),\mathbf{c}_p\big]^{(-1)^{(p+1)}},$$
an element of $\mathbb{C}^{*}/\{\pm 1\}$.
\end{definition}

The twisted Reidemeister torsion depends on the choice of homology bases $\mathbf{h}_p$; however, it is independent of the conjugacy class of $\rho$, the choice of lifts $\widetilde{e}^p_j$, and the choice of the $\mathcal{B}$-orthonormal basis $\mathcal{A}$ (see \cite[Section 0.2, Remarks a1, a2]{Porti}), and it does not depend on the auxiliary bases $\mathbf{b}_p$ or sections $s_p,\ell_p$ (see \cite[Section 3, p.\,365]{Milnor}).

By Milnor's subdivision-invariance result \cite{Milnor}, the twisted Reidemeister torsion defines an invariant of the underlying manifold, independent of triangulation.

\begin{definition}
Let $M$ be a smooth compact $n$-manifold ($n=1,2,3$) with a triangulation $K$, let $\rho:\pi_1(M)\to G$ be a representation, and let $\{\mathbf{h}_p\}_{p=0}^{n}$ be homology bases. The twisted Reidemeister torsion of $M$ is
$$\mathbb{T}_{\rho}(M,\{\mathbf{h}_p\}_{p=0}^{n})=\mathbb{T}(C_{\ast}(K;\mathfrak{g}_{\mathrm{Ad}_\rho}),\{\mathbf{c}_p\}_{p=0}^{n},\{\mathbf{h}_p\}_{p=0}^{n}),$$
which is independent of the triangulation $K$ of $M$.
\end{definition}

The next result records the multiplicativity of twisted Reidemeister torsion under gluing; it is the key computational tool used throughout the paper.

\begin{theorem}[See \cite{Porti}]\label{prt1}
Let $X$ be a compact CW-complex with subcomplexes $X_1,X_2\subset X$ such that $X=X_1\cup X_2$ and $Y=X_1\cap X_2$, and let $Y_1,\ldots,Y_k$ be the connected components of $Y$. For the inclusions
$$Y\overset{i_\nu}{\hookrightarrow} X_\nu\overset{j_\nu}{\hookrightarrow} X,$$
let $\rho:\pi_1(X)\to G$ be a representation, with restrictions $\rho_{|_{X_\nu}}:\pi_1(X_\nu)\to G$ and $\rho_{|_{Y_\mu}}:\pi_1(Y_\mu)\to G$, $\nu=1,2$, $\mu=1,\ldots,k$. Then:

\begin{itemize}
\item[(i)] The sequence
\begin{eqnarray*}
&&0\to\bigoplus_{\mu}C_{\ast}(Y_\mu;\mathfrak{g}_{\mathrm{Ad}_{\rho_{|_{Y_\mu}}}})
\overset{(i_1)_\#\oplus(i_2)_\#}{\longrightarrow}
C_{\ast}(X_1;\mathfrak{g}_{\mathrm{Ad}_{\rho_{|_{X_1}}}})\oplus C_{\ast}(X_2;\mathfrak{g}_{\mathrm{Ad}_{\rho_{|_{X_2}}}})\\
&&\qquad\qquad\qquad\overset{(j_1)_\#-(j_2)_\#}{\longrightarrow}
C_{\ast}(X;\mathfrak{g}_{\mathrm{Ad}_\rho})\to 0
\end{eqnarray*}
is short exact.

\item[(ii)] There is a Mayer--Vietoris long exact sequence in twisted homology,
\begin{eqnarray}\label{es12es}
\mathcal{H}_{\ast}:\ \cdots&\to&\bigoplus_{\mu}H_i(Y_\mu;\mathfrak{g}_{\mathrm{Ad}_{\rho_{|_{Y_\mu}}}})
\overset{i_1^\ast\oplus i_2^\ast}{\longrightarrow}
H_i(X_1;\mathfrak{g}_{\mathrm{Ad}_{\rho_{|_{X_1}}}})\oplus H_i(X_2;\mathfrak{g}_{\mathrm{Ad}_{\rho_{|_{X_2}}}}) \nonumber\\
&\overset{j_1^\ast-j_2^\ast}{\longrightarrow}&
H_i(X;\mathfrak{g}_{\mathrm{Ad}_\rho})\to\bigoplus_{\mu}H_{i-1}(Y_\mu;\mathfrak{g}_{\mathrm{Ad}_{\rho_{|_{Y_\mu}}}})\to\cdots
\end{eqnarray}

\item[(iii)] For chosen homology bases $\mathbf{h}^X_i$, $\mathbf{h}_i^{X_1}$, $\mathbf{h}^{X_2}_i$, $\mathbf{h}^{Y_\mu}_i$, viewing \eqref{es12es} as a based chain complex,
\begin{eqnarray*}
\mathbb{T}_{\rho_{|_{X_1}}}(X_1,\{\mathbf{h}^{X_1}_i\})\;\mathbb{T}_{\rho_{|_{X_2}}}(X_2,\{\mathbf{h}^{X_2}_i\})
&=&\mathbb{T}_{\rho}(X,\{\mathbf{h}^X_i\})\prod_{\mu=1}^k\mathbb{T}_{\rho_{|_{Y_\mu}}}(Y_\mu,\{\mathbf{h}^{Y_\mu}_i\})\\
&&\times\ \mathbb{T}(\mathcal{H}_{\ast},\{\mathbf{h}_{\ast\ast}\}).
\end{eqnarray*}
\end{itemize}
\end{theorem}

We refer to \cite[Theorem 3.2]{Milnor} and \cite[Section 0.4, Proposition 0.11]{Porti} for the proof of Theorem \ref{prt1}.

The factor $\mathbb{T}(\mathcal{H}_{\ast},\{\mathbf{h}_{\ast\ast}\})$ in Theorem \ref{prt1}(iii) is called the \emph{corrective term}; it can be computed as follows.

\begin{lemma}[See \cite{EDErdal}]\label{correctivetermlemma}
Let $X$ be an $n$-dimensional CW-complex. Then the corrective term satisfies
\begin{equation*}
\mathbb{T}\left(\mathcal{H}_{\ast},\{\mathbf{h}_p\}_{p=0}^{3n+2},\{0\}_{p=0}^{3n+2}\right)
=\prod_{p=0}^{3n+2}\big[\mathbf{h}'_p,\mathbf{h}_p\big]^{(-1)^{(p+1)}},
\end{equation*}
where $\mathbf{h}'_p=\mathbf{b}_p\sqcup s_p(\mathbf{b}_{p-1})$.
\end{lemma}

\subsection{The Twisted Riedemeister Torsion of Disk Sum of Manifolds}
Let $M$ be a $3$-manifold with connected non-vacuous boundary which is not spherical, and let $N$ be another $3$-manifold with vacuous boundary. Recall the \emph{disk sum} $M\bigtriangleup N^{\ast}$, where $N^{\ast}$ is $N$ with the interior of a $3$-cell removed: it is formed by pasting a $2$-cell on the boundary of $M$ to a $2$-cell on the boundary of $N^{\ast}$. The operation $\bigtriangleup$ is well-defined, associative and commutative up to homeomorphism, and by \cite{JLGROSS} every compact $3$-manifold admits a unique decomposition into $\bigtriangleup$-prime pieces. Using this decomposition, the following theorem shows that twisted Reidemeister torsion is multiplicative under the disk sum operation.

\begin{theorem}[See Theorem 1.2 in \cite{EDErdal}]\label{theo1}
Consider the disk sum $M=\bigtriangleup_{i=1}^{n}M_i$, and denote by $\mathbb{D}_j^2$ ($j=1,\ldots,n-1$) the disks on the boundaries of $M_j$ and $M_{j+1}$ identified in the construction of $M$, assumed pairwise disjoint in their interiors. Let $\varrho:\pi_1(M)\to G$ be a representation with restrictions $\varrho_{|_{\pi_1(M_i)}}=\psi_i$. For a given basis $\mathbf{h}^{M}_p$ of $H_p(M;\mathfrak{g}_{\mathrm{Ad}_\varrho})$, and the basis $\mathbf{h}_0^{\mathbb{D}_j^2}=f^j_\ast(\varphi_0(\mathbf{c}_0))$ of $H_0(\mathbb{D}_j^2;\mathfrak{g}_{\mathrm{Ad}_{\varrho_{|_{\mathbb{D}_j^2}}}})$, there exist bases $\mathbf{h}^{M_i}_p$ of $H_p(M_i;\mathfrak{g}_{\mathrm{Ad}_{\psi_i}})$ for each $i$ such that
\begin{eqnarray*}
\mathbb{T}_{\varrho}(M,\{\mathbf{h}^{M}_p\}_{p=0}^3)=\prod_{i=1}^{n}\mathbb{T}_{\psi_i}(M_i,\{\mathbf{h}_p^{M_i}\}_{p=0}^{3}),
\end{eqnarray*}
where $f^j_\ast$ is induced by the simple homotopy equivalence $f^j:\{pt.\}\to\mathbb{D}_j^2$, $\varphi_0:Z_0(\{pt.\};\mathrm{Ad}_\varphi)\to H_0(\{pt.\};\mathrm{Ad}_\varphi)$ is the natural projection, and $\mathbf{c}_0^j$ is the geometric basis of $C_0(\{pt.\};\mathrm{Ad}_\varphi)$.
\end{theorem}

\section{ The Twisted Riedemeister Torsion of Connected Sum of Manifolds}
\label{Sec:3}

Combining Theorem \ref{theo1} with the homeomorphism $M\# N\approx M\bigtriangleup N^{\ast}$ of \cite{JLGROSS}, we now prove our main result.\\

\noindent \textbf{Proof of Theorem \ref{primetheo}.} Let $\varphi:M\# N\to M\bigtriangleup N^{\ast}$ be the homeomorphism from \cite{JLGROSS}. At the level of homology this induces isomorphisms
$$[\varphi_i]:H_i(M\#N)\longrightarrow H_i(M\bigtriangleup N^{\ast}),\qquad i=0,1,2,3,$$
and, via the induced isomorphism $[\varphi^{-1}]:\pi_1(M\bigtriangleup N^{\ast})\to\pi_1(M\#N)$, a representation
$$\sigma=[\varphi^{-1}]\circ\varrho:\pi_1(M\bigtriangleup N^{\ast})\longrightarrow G.$$

Let $\{\mathbf{h}_i^{M\#N}\}_{i=0}^3$ be a given homology basis of 
$H_i(M\# N;\mathfrak{g}_{\mathrm{Ad}_{\varrho}})$. Since Reidemeister torsion is a topological invariant, the homeomorphism $\varphi$ yields
\begin{equation}\label{prmequatn4}
\mathbb{T}_\varrho\big(M\#N,\{\mathbf{h}_i^{M\#N}\}_{i=0}^3\big)
=\mathbb{T}_\sigma\big(M\bigtriangleup N^{\ast},\{[\varphi_i](\mathbf{h}_i^{M\#N})\}_{i=0}^3\big).
\end{equation}

The disk sum $M\bigtriangleup N^{\ast}$ is obtained by gluing a topological disk $\mathbb{D}^2$ on $\partial M$ to a topological disk $\mathbb{D}^2$ on $\partial N^{\ast}$. By Theorem \ref{prt1}, computing the twisted torsion of $M\bigtriangleup N^{\ast}$ via this gluing involves the factor $\mathbb{T}_{\sigma_{|_{\mathbb{D}^2}}}(\mathbb{D}^2,\{\mathbf{h}_0^{\mathbb{D}^2}\})$. Choosing $\mathbf{h}_0^{\mathbb{D}^2}=f_\ast(\varphi_0(\mathbf{c}_0))$ as in Theorem \ref{theo1}, where $f:\{pt.\}\to\mathbb{D}^2$ is the simple homotopy equivalence, Lemma \cite[Lemma 3.1]{EDErdal} gives
$$\mathbb{T}_{\sigma_{|_{\mathbb{D}^2}}}\big(\mathbb{D}^2,\{\mathbf{h}_0^{\mathbb{D}^2}\}\big)=1.$$

Let $\psi_M=\sigma_{|_{\pi_1(M)}}$ and $\psi_{N^{\ast}}=\sigma_{|_{\pi_1(N^{\ast})}}$ be the restrictions of $\sigma$. Applying Theorem \ref{theo1} to the disk sum $M\bigtriangleup N^{\ast}$ with basis $\{[\varphi_i](\mathbf{h}_i^{M\#N})\}_{i=0}^3$, there exist homology bases $\{\mathbf{h}^{M}_i\}_{i=0}^3$ and $\{\mathbf{h}^{N^{\ast}}_i\}_{i=0}^3$ of $H_i(M;\mathfrak{g}_{\mathrm{Ad}_{\psi_M}})$ and $H_i(N^{\ast};\mathfrak{g}_{\mathrm{Ad}_{\psi_{N^{\ast}}}})$ such that
\begin{equation}\label{eqnlst}
\mathbb{T}_\sigma\big(M\bigtriangleup N^{\ast},\{[\varphi_i](\mathbf{h}_i^{M\#N})\}_{i=0}^3\big)
=\mathbb{T}_{\psi_M}\big(M,\{\mathbf{h}_i^{M}\}_{i=0}^3\big)\;\mathbb{T}_{\psi_{N^{\ast}}}\big(N^{\ast},\{\mathbf{h}_i^{N^{\ast}}\}_{i=0}^3\big).
\end{equation}

Combining \eqref{prmequatn4} and \eqref{eqnlst} gives precisely formula \eqref{mainformula}, and shows that twisted Reidemeister torsion is multiplicative under connected sum, \emph{without any corrective term}:
\begin{eqnarray*}
\mathbb{T}_\varrho\big(M\#N,\{\mathbf{h}_i^{M\#N}\}_{i=0}^3\big)
=\mathbb{T}_{\psi_M}\big(M,\{\mathbf{h}_i^{M}\}_{i=0}^3\big)\;\mathbb{T}_{\psi_{N^{\ast}}}\big(N^{\ast},\{\mathbf{h}_i^{N^{\ast}}\}_{i=0}^3\big).
\end{eqnarray*}
This completes the proof of Theorem \ref{primetheo}. $\hfill\blacksquare$

\end{document}